\documentclass[a4paper,11pt]{amsart}

\usepackage{rotating,pdflscape,color,amssymb,hyperref,subfigure,psfrag,amsmath,eufrak,bbm,epsfig}
\author[Florent Benaych-Georges]{Florent Benaych-Georges}\address{Florent Benaych-Georges, MAP 5, UMR CNRS 8145 - Universit\'e Paris Descartes, 45 rue des Saints-P\`eres 75270 Paris cedex~6, France and CMAP, \'Ecole Polytechnique, route de Saclay, 91128 Palaiseau Cedex, France. URL: {\rm  http://www.cmapx.polytechnique.fr/$\sim$benaych/}} \email{florent.benaych@gmail.com}
\title[Eigenvectors of Wigner matrices: universality of global fluctuations]{A universality result for the global fluctuations of the eigenvectors of Wigner matrices}
\date{\today}
\thanks{This work was partially supported by the \emph{Agence Nationale de la Recherche} grant ANR-08-BLAN-0311-03 and   partly accomplished during the author's stay at  New York University Abu Dhabi, Abu Dhabi, U.A.E.}

\newcommand{\op}{\operatorname}
\newcommand{\re}{\ref}
\newcommand{\la}{\label}
\newcommand{\dist}{\operatorname{distribution}}

\newcommand{\sh}{\#}

\newcommand{\si}{\sigma}

\newcommand{\scl}{\operatorname{semicircle}}

\newcommand{\Co}{\operatorname{Cov}}

\newcommand{\diag}{\operatorname{diag}}

\newcommand{\Part}{\operatorname{Part}}

\newcommand{\Pro}{\mathbb{P}}

\newcommand{\Tr}{\operatorname{Tr}}

\newcommand{\ninf}{\underset{n\to\infty}{\longrightarrow}}

\newcommand{\ssi}{if and only if }

\newcommand{\ensn}{\{1,\ldots,n\}}
\newcommand{\E}{\mathbb{E}}

\newcommand{\R}{\mathbb{R}}
\newcommand{\C}{\mathbb{C}}

\newcommand{\ud}{\mathrm{d}}

\newcommand{\pro}{probability }

\newcommand{\f}{\frac}
\newcommand{\ff}{\frac{1}}
\newcommand{\lf}{\left}
\newcommand{\ri}{\right}

\newcommand{\st}{such that }
\newcommand{\stm}{\textrm{s.t. }}

\newcommand{\lam}{\lambda}

\newcommand{\ste}{\, ;\, }
\newcommand{\mc}{\mathcal }

\newcommand{\eps}{\varepsilon}

\newcommand{\al}{\alpha}

\newcommand{\ovl}{\overline}

\newcommand{\bbm}{\begin{bmatrix}}
\newcommand{\ebm}{\end{bmatrix}}
\newcommand{\bes}{\begin{equation*}}
\newcommand{\ees}{\end{equation*}}
\newcommand{\be}{\begin{equation}}
\newcommand{\ee}{\end{equation}}
\newcommand{\beqy}{\begin{eqnarray}}
\newcommand{\eeqy}{\end{eqnarray}}
\newcommand{\beq}{\begin{eqnarray*}}
\newcommand{\eeq}{\end{eqnarray*}}
\newcommand{\one}{\mathbbm{1}}
\newcommand{\lto}{\longrightarrow}

\newtheorem{Th}{Theorem}[section]

\newtheorem{assum}[Th]{Assumption} 
\newtheorem{propo}[Th]{Proposition} 

\newtheorem{lem}[Th]{Lemma}

\newtheorem{rmq}[Th]{Remark}

\newenvironment{pr}{\noindent {\bf Proof. }}{\hfill $\square$}

\long\def\symbolfootnote[#1]#2{\begingroup
\def\thefootnote{\fnsymbol{footnote}}\footnote[#1]{#2}\endgroup}

\begin{document}
\maketitle

\begin{abstract}We prove that for $[u_{i,j}]_{i,j=1}^n$  the eigenvectors matrix of a Wigner matrix, under some moments conditions, the bivariate random process $$\lf(\sum_{\substack{1\le i\le ns,\\ 1\le j\le nt}}(|u_{i,j}|^2-1/n)\ri)_{(s,t)\in [0,1]^2}$$
converges in distribution to a bivariate Brownian bridge. This result has already been proved for GOE and GUE matrices. It is conjectured here that the necessary and sufficient condition, for the result to be true for a general Wigner matrix, is the matching of the moments of orders $1$, $2$ and $4$ of the entries of the Wigner with the ones of a GOE or GUE matrix.  Surprisingly, the third moment of the entries of the Wigner matrix has no influence on the limit distribution.\end{abstract}


\section{Introduction}
It is well known that the matrix $U_n=[u_{i,j}]_{i,j=1}^n$ whose columns are the eigenvectors of a GOE or GUE matrix $X_n$ can be chosen to be distributed according to the Haar measure on the orthogonal or unitary group. As a consequence, much can be said about the $u_{i,j}$'s: their joint moments can be computed via the so-called {\it Weingarten calculus} developed in \cite{collinsIMRN,collins-sniady06}, any finite (or not too large) set of $u_{i,j}$'s can be approximated, as $n\to\infty$, by independent Gaussian variables (see \cite{jiang06,meckessourav08}) and the global asymptotic fluctuations of the $|u_{i,j}|$'s are governed by a theorem of Donati-Martin and Rouault, who proved in \cite{cat-alain10}  that as $n\to\infty$, the bivariate c\`adl\`ag process 
   $$\lf(B^n_{s,t}:=\sqrt{\f{\beta}{2}}\sum_{\substack{1\le i\le ns,\\ 1\le j\le nt}}(|u_{i,j}|^2-1/n)\ri)_{(s,t)\in [0,1]^2}$$(where    $\beta=1$ in the real case and $\beta=2$ in the complex case) converges in distribution, for the Skorokhod topology,  to the {\it bivariate Brownian bridge}, i.e. the centered continuous Gaussian process $(B_{s,t})_{(s,t)\in [0,1]^2}$  with covariance \be\label{cov_bb}\E[B_{s,t}B_{s',t'}]=(\min\{s,s'\} -ss')(\min \{t,t'\}-tt').\ee  
 
A natural question is the following:

 \begin{quote}What can be said beyond the Gaussian case, when the entries of the Wigner matrix $X_n$ are general random variables ?\end{quote}
 
 For a general Wigner matrix\footnote{A  {\it Wigner matrix} is a real symmetric or Hermitian random matrix with independent, centered entries whose variance is one. Its {\it atom distributions} are the distributions of its entries.}, the exact distribution of the matrix $U_n$ cannot be computed and few works had been devoted to this subject until quite recently. One of the reasons is that while the eigenvalues of an Hermitian  matrix 	admit variational characterizations as {\it extremums} of certain functions, the eigenvectors can be characterized as the {\it argmax}  of these functions, hence are more sensitive to perturbations of the entries of the matrix. 
 However, in the last three years, the eigenvectors of general Wigner matrices have been the object of a growing interest, due in part to some relations with the universality conjecture for the eigenvalues. In several papers (see, among others, \cite{ErdosSchleinYau1,ErdosSchleinYau2,ErdosSchleinYau3}), a {\it delocalization property} was shown for the eigenvectors of random matrices. More recently, 
 Knowles and Yin in \cite{KnowlesYinEig} and  Tao and Vu in \cite{Tao-Vu_1103.2801v1}  proved that if the first four moments of the atom distributions of $X_n$ coincide with the ones of a GOE or GUE matrix, then under some tail assumptions on these distributions,    the $u_{i,j}$'s can be approximated by independent Gaussian   variables as long as we only consider  a finite  (or not too large) set of $u_{i,j}$'s. 
 
 In this paper, we consider the \emph{global behavior}  of the $|u_{i,j}|$'s. By \emph{global behavior}, we mean that we consider functionals of the  $|u_{i,j}|$'s that involve all of them, at the difference of the works  of Knowles and Yin in \cite{KnowlesYinEig} and  Tao and Vu in \cite{Tao-Vu_1103.2801v1}. Our work is in the same vein as Silverstein's paper \cite{silverstein-AOP-eigenvectors} or Bai and Pan's paper \cite{baipanJStatPhys2012}. We prove (Theorem \ref{main_th0}) that for   Wigner matrices whose entries have moments of all orders, the process $(B^n_{s,t})_{(s,t)\in [0,1]^2}$ has a limit in a weaker sense than for the Skorokhod topology and that this weak limit is the bivariate Brownian bridge \ssi the off-diagonal entries of the matrix have the same fourth moment as the GOE or GUE matrix (quite surprisingly, no hypothesis on the third moment is necessary). Under some additional hypotheses on the atom distributions (more coinciding moments and   continuity), we prove the convergence for the Skorokhod topology (Theorem \ref{main_th}).
 
 This result was conjectured by Djalil Chafa\"\i, who also conjectures the same kind of universality for  unitary matrices appearing in other standard decompositions, such as the singular values decomposition  or the Housholder decomposition of   matrices with no symmetry, as long as the matrix considered has i.i.d. entries with    first moments agreeing with the ones of Gaussian variables.  It would also be interesting to consider the same type of question in the context of {\it band matrices}, connecting this problem  with the so-called {\it Anderson conjecture} (see \emph{e.g.} the works of Erd\"os and Knowles \cite{erdosknowlesband1, erdosknowlesband2}, of Schenker \cite{schenker_band} or of Sodin \cite{s_sodin_band_matrices}, or, for a short introduction,  the blog note by   Chafa\"{\i} \cite{chafai-anderson}).  
 
 The paper is organized as follows. The main results are stated in Section \ref{sec:beyond}, where we also make some comments on their hypotheses ;  outlines of the proofs and the formal proofs are given in Section \ref{secproof} ; and Section \re{sectechresults} is devoted to the definitions of the functional spaces and their topologies and to the proofs of several technical results needed in Section \ref{secproof}.
\\

\noindent{\bf Acknowledgements: }It is a pleasure for the author to thank Djalil Chafa{\"\i}  for having generously  pointed out   this  problem  to him  and for the nice discussions we had about it. We also would like to thank Alice Guionnet, for her availability and her precious advices, and Terry Tao, who kindly and patiently answered  several naive questions asked by the author on his blog. At last, we would like to thank the anonymous referee who pointed a mistake in the proof of Proposition \ref{bound8411} and brought the paper \cite{baipanJStatPhys2012} to the attention of the author.
 
\section{Main results}\label{sec:beyond}

For each $n$, let us consider a  real symmetric or Hermitian random matrix $$X_{n}:= \ff{\sqrt{n}}[x_{i,j}^{(n)} ]_{i,j=1}^n.$$ For notational brevity,  $x_{i,j}^{(n)}$ will be denoted by $x_{i,j}$.

 Let us denote by $\lam_1\le\cdots\le \lam_n$ the eigenvalues of $X_n$  and consider an orthogonal or unitary   matrix $U_n=[u_{i,j}]_{i,j=1}^n$ \st $$X_n=U_n\diag(\lam_1, \ldots, \lam_n) U_n^*.$$ Note that   $U_n$ is not uniquely defined. However, one can choose it in any measurable way.

We define the bivariate c\`adl\`ag process
 $$\lf(B^n_{s,t}:=\sqrt{\f{\beta}{2}}\sum_{\substack{1\le i\le ns,\\ 1\le j\le nt}}(|u_{i,j}|^2-1/n)\ri)_{(s,t)\in [0,1]^2},$$where $\beta=1$ in the real case and  $\beta=2$ in the complex case.  
 
  \begin{assum}\label{hyp_vector:indep_entries} For each $n$, the random variables $x_{i,j}$'s are independent (up to the symmetry), with the same distribution on the diagonal and the same distribution above the diagonal.\end{assum}
 
    \begin{assum}\label{hyp_vector:moments2} For each $k\ge 1$, $\sup_n \E[|x_{1,1}|^{k}+|x_{1,2}|^{k}]<\infty$. 
  Moreover,     \be\label{bound_moments31311pp2} \E[x_{1,1}]= 
  \E[x_{1,2}]=0,\qquad\qquad
  \E[|x_{1,2}|^2]=1\ee 
  and  $\E[|x_{1,2}|^4]$ has a limit  as $n\to\infty$.
  \end{assum}

 The  bivariate Brownian bridge has been defined in the introduction  and the definitions of the functional spaces and their topologies can be found in Section \ref{sec:espaces+topo}. 
 
    \begin{Th}\label{main_th0} Suppose that Assumptions \ref{hyp_vector:indep_entries} and \ref{hyp_vector:moments2}  are satisfied. 
 Then the sequence $$(\operatorname{distribution}(B^n))_{n\ge 1}$$  has a unique possible accumulation  point supported by $C([0,1]^2)$. This accumulation  point  is the distribution of a centered Gaussian process which depends on the distributions of the $x_{i,j}$'s only through  $\lim_{n\to\infty}\E[|x_{1,2}|^4]$, and which is the  bivariate Brownian bridge  \ssi $\lim_{n\to\infty} \E[|x_{1,2}|^4]=4-\beta$.  
 \end{Th}
 
 More precisions about the way the unique possible accumulation point depends on the   fourth moment of the entries are given in Remark \re{prec_lim_2611}. 
 
 To get a stronger statement where the convergence in distribution to the bivariate Brownian bridge is actually stated, one needs stronger hypotheses.

    \begin{assum}\label{hyp:density} The distributions of the entries of $X_n$ are absolutely continuous with respect to the Lebesgue measure.\end{assum} 
 The following hypothesis depends on an integer $m\ge 2$.
  \begin{assum}\label{hyp_vector:moments} For each $k\ge 1$, $\sup_n \E[|x_{1,1}|^{k}+|x_{1,2}|^{k}]<\infty$. 
  Moreover,   for each $r,s\ge 0$,  \be\label{bound_moments31311pp}r+s\le m-2\;\Longrightarrow\; \E[\Re(x_{1,1})^r\Im(x_{1,1})^s]=\E[\Re(g_{1,1})^r\Im(g_{1,1})^s] \ee and \be\label{bound_moments31311pq}r+s\le m\;\Longrightarrow \;\E[\Re(x_{1,2})^r\Im(x_{1,2})^s]=\E[\Re(g_{1,2})^r\Im(g_{1,2})^s].\ee
  where the $g_{i,j}$'s are the entries of a standard GOE or GUE matrix.\end{assum}

 \begin{Th}\label{main_th} Suppose that Assumptions \ref{hyp_vector:indep_entries} and   \ref{hyp:density} are satisfied, as well as Assumption    \ref{hyp_vector:moments} for $m=12$. 
 Then, as $n\to\infty$, the bivariate process $B^n$ converges in distribution,  for the Skorokhod topology in $D([0,1]^2)$,  to the bivariate Brownian bridge. 
 \end{Th}

  \begin{rmq}{\rm {\bf Weakening of the assumptions.}  As explained above, the distributions of the entries $x_{i,j}=x_{i,j}^{(n)}$ are allowed to depend on $n$.
  In order to remove Assumption \ref{hyp:density} below (which we did not manage to do yet), it mights be useful to weaken Assumptions \ref{hyp_vector:moments2} and \ref{hyp_vector:moments} in the following way:   one can easily see   that the proofs of Theorems \ref{main_th0} and \ref{main_th} and Proposition \ref{bound8411} still work if, in Equations \eqref{bound_moments31311pp2},  
   \eqref{bound_moments31311pp} and \eqref{bound_moments31311pq},  one   replaces the identities by the same identities up to an error which is $O(n^{-\al})$ for all $\al>0$.}\end{rmq}

  \begin{rmq}\la{prec_lim_2611}{\rm {\bf Complements on Theorem \re{main_th0}.}  One can wonder how the unique accumulation point mentioned in Theorem \re{main_th0} depends on the fourth moment of the entries of  $X_n$. Let  $G:=(G_{s,t})_{(s,t)\in [0,1]^2}$ be distributed according to this distribution. We know that $(G_{s,t})_{(s,t)\in [0,1]^2}$ is the bivariate Brownian bridge only in the case where $\lim_{n\to\infty}\E[|x_{1,2}|^4]=4-\beta$. In the other cases, defining $F_{\scl}$ as the cumulative distribution function of the semicircle law, the covariance of the centered Gaussian process \be\la{moments_pro_sto_lim}\lf(\int_{u=-2}^2u^kG_{s,F_{\scl}(u)}\ud u\ri)_{s\in [0,1], k\ge 0}\ee can be computed thanks to Lemma \re{moments-mesuredemassenulle} and Proposition \re{conv_moments}. By Lemma \re{moments_pro_sto}, it determines completely the distribution of the process $G$. However, making the covariance of $G$ explicit out of the covariance of the process of \eqref{moments_pro_sto_lim} is  a very delicate problem, and we shall only stay at a quite vague level, saying that it can be deduced  from the proof of Proposition \ref{conv_moments} that the variances of the one-dimensional marginals of $G$ are increasing functions of $\lim_{n\to\infty}\E[|x_{1,2}|^4]$. More specifically, one can deduce from Lemma \re{moments-mesuredemassenulle} and Proposition \ref{conv_moments} that for all $0\le s_1,s_2\le 1$, $$\op{Cov}\lf(\int_{u=-2}^2\!\!\!\!u^2G_{s_1,F_{\scl}(u)}\ud u, \int_{u=-2}^2\!\!\!\!u^2G_{s_2,F_{\scl}(u)}\ud u\ri)=\f{\E[|x_{1,2}|^4]-1}{4}(\min\{s_1,s_2\}-s_1s_2).$$}\end{rmq}\vspace{0.3cm}
 
   \begin{rmq}{\rm {\bf Comments on the hypotheses of Theorem \ref{main_th} (1).} In order to prove the convergence in the Skorokhod topology, we had to make several hypotheses on the atom distributions: absolute   continuity, moments of all orders  and   coincidence of their 10 (on the diagonal) and 12 (above the diagonal) first moments   with the ones of a GOE or GUE matrix. We needed these assumptions    to control the discontinuities of the process $B^n$. Even though these hypotheses might not be optimal (especially the continuity one), a bound on the tails of the atom distributions seems to be necessary to avoid   too large variations of the process $B^n$. Indeed, as illustrated  by Figure \ref{Fig:eigenvectors}, 
    \begin{figure}[h!]
\begin{center}
\includegraphics[width=15cm,height=7cm]{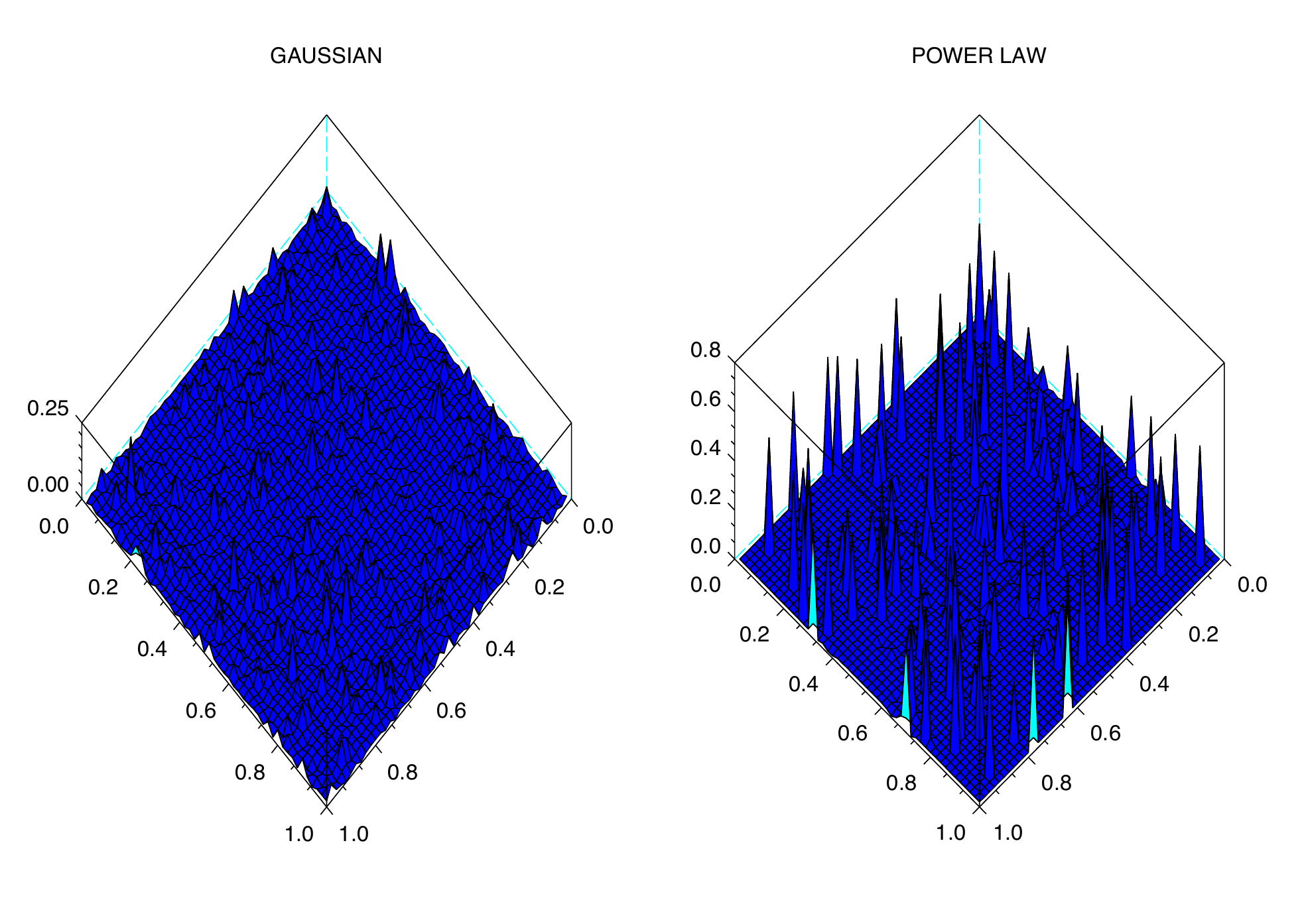}
\caption{{\bf Influence of the tails of the atom distributions of $X_n$ on the $|u_{i,j}|$'s:} Plot of the map $(i/n,j/n)\longmapsto ||u_{i,j}|^2-1/n|$ for two different choices of atom distributions. {\it Left:} GOE matrix.  {\it Right:}  Wigner matrix with   atom distribution admitting moments only up to order $2+\eps$ for a small $\eps$.  For both pictures, the matrices are $n\times n$ with $n=50$.}\label{Fig:eigenvectors}
\end{center}
\end{figure}
     for a GOE matrix  ({\it left} picture), $|u_{i,j}|^2$ is close to $1/n$ for all $i,j$ with high probability,  whereas when the atom distributions have not more than a second moment ({\it right} picture), the matrix $X_n$ looks more like a {\it sparse matrix}, and so does $U_n$, which implies that for certain $(i,j)$'s, $|u_{i,j}|^2-1/n$ is not small enough. Since $|u_{i,j}|^2-1/n$ is the   jump of the process $B^n$ at $(s,t)=(i/n,j/n)$, this could be an obstruction to the existence of a continuous limit for the process $B^n$. That being said, we have hopes to prove the theorem under Assumption \ref{hyp_vector:moments} for $m=4$ instead of $m=12$ (see Remark \ref{22041112h30} bellow).}\end{rmq}

Note that it follows from the previous theorem that for all $0\le s<s'\le 1$ and $0\le t<t'\le 1$, the sequence of random variables $$\ff{\sqrt{(s'-s)(t'-t)}}\sum_{\substack{ns<i\le ns'\\ nt<j\le nt'}}(|u_{i,j}|^2-1/n)$$ admits a limit in distribution as $n\to\infty$, hence is bounded in \pro (in the sense of \cite[Def. 1.1]{taovucircularlaw}). In the same way, it follows from \cite{KnowlesYinEig} and    \cite{Tao-Vu_1103.2801v1} that the sequence $n|u_{i,j}|^2-1$ is bounded in probability. In the next proposition, we improve these assertions by making them uniform on $s,s',t,t',i,j$ and upgrading them to the $L^2$ and $L^4$ levels. The proof of the proposition is postponed to Section \re{prooflem641114h}.

 \begin{propo}\label{bound8411} Suppose that Assumptions \ref{hyp_vector:indep_entries},   \ref{hyp:density} and  \ref{hyp_vector:moments} for $m=4$  (resp. $m=8$) are satisfied. 
 Then as $n\to\infty$, the sequence   \be\label{terme_maj_22042011} n|u_{i,j}|^2-1\qquad\textrm{ (resp.  }\qquad\ff{\sqrt{(s'-s)(t'-t)}}\sum_{\substack{ns<i\le ns'\\ nt<j\le nt'}}(|u_{i,j}|^2-1/n)\textrm{ )}\ee is bounded for the  $L^4$ (resp. $L^2$) norm, uniformly in $i,j$ (resp. $s<s',t<t',$).
 \end{propo}
 
  \begin{rmq}{\rm {\bf Comments on the hypotheses of Theorem \ref{main_th} (2).}\label{22041112h30} This proposition  is almost sufficient to apply the tightness criterion that we use in this paper. Would  the second term of  \eqref{terme_maj_22042011} have been bounded for the $L^{2+\eps}$ norm (instead of $L^2$), Assumption \ref{hyp_vector:moments} for $m=8$ would have been enough to prove that $B^n$ converges in distribution,  for the Skorokhod topology in $D([0,1]^2)$,  to the bivariate Brownian bridge.}\end{rmq}

\section{Proofs of Theorems \ref{main_th0} and \ref{main_th}}\la{secproof}

\subsection{Outline of the proofs}
 Firstly, Theorem \ref{main_th} can be deduced from Theorem \ref{main_th0} by proving that the sequence $(distribution(B^n))_{n\ge 1}$ is tight and  only has $C([0,1]^2)$-supported accumulation points.   This can be done via   some upper bounds on the fourth  moment of the increments of $B^n$ and on its jumps (or discontinuities). These bounds are given   in the proof of  Lemma \ref{tightness_lemma} below, and rely  on the existing bounds in the case where $X_n$ is a GOE or GUE matrix and on the ``one-by-one entries replacement method" developed by Terence Tao and Van Vu in recent papers, such as \cite{Tao-Vu_0906.0510,Tao-Vu_1103.2801v1}.

Secondly, the proof of Theorem \ref{main_th0} relies on the following remark, inspired by some ideas of Jack Silverstein (see \cite[Chap. 10]{bai-silver-book} and \cite{silverstein-AOP-eigenvectors}): even though we do not have any ``direct access" to the eigenvectors of $X_n$,  we have access to the process $(B^n_{s, F_{\mu_{X_n}}(u)})_{s\in [0,1], u\in \R}$, for  $F_{\mu_{X_n}}(u):=\ff{n}\sh\{i\ste \lam_i\le u\}$. 
Indeed,   $$
B^n_{s, F_{\mu_{X_n}}(u)}=\sqrt{\f{\beta}{2}}\sum_{1\le i\le ns}\sum_{\substack{1\le j\le n\\ \stm \lam_j\le u}}(|u_{i,j}|^2-1/n),$$hence for all fixed $s\in [0,1]$, the function $u\in \R \longmapsto B^n_{s, F_{\mu_{X_n}}(u)}$ is the cumulative distribution function of the  signed measure\footnote{Note that the random \pro measures $\sum_{j=1}^n|u_{i,j}|^2\delta_{\lam_j}$ have been studied in \cite{baipanJStatPhys2012} (convergence to the semicircle law and fluctuations).  Here, we are somehow considering   the fluctuations of large sums of such random \pro measures ($1\leq i\leq ns$) around the empirical spectral law $\ff{n}\sum_{j=1}^n \delta_{\lam_j}$.} 
 \be\label{441118hintro} \sqrt{\f{\beta}{2}}\sum_{1\le i\le ns}\sum_{j=1}^n(|u_{i,j}|^2-1/n)\delta_{\lam_j},\ee which  can   be studied via its moments $$\sum_{1\le i\le ns}\lf(e_i^*X_n^{k}e_i-\ff{n}\Tr X_n^{k}\ri)\qquad \textrm{ ($k\ge 1$)},$$   the $e_i$'s being the vectors of the canonical basis. From the asymptotic behavior of the {\it moments} of the 
 signed measure 
 of \eqref{441118hintro}, one can then find out the 
  asymptotic behavior of its {\it cumulative distribution function}.  
 
 Once the asymptotic distribution  of the process $(B^n_{s, F_{\mu_{X_n}}(u)})_{s\in [0,1], u\in \R}$ identified, one can obtain the asymptotic distribution of the process $(B^n_{s, t})_{s\in [0,1], t\in[0,1]}$ because the function $F_{\mu_{X_n}}$ tends to the (non random) cumulative distribution function $F_{\scl}$ of the semicircle law. 

\subsection{Formal proofs}

By a standard {\it tightness + uniqueness of the accumulation point} argument, Theorem  \ref{main_th} will   follow from the following lemma and  Theorem  \ref{main_th0}. The proof of the lemma is postponed to Section \re{prooflem641114h}. 

\begin{lem}\label{tightness_lemma}Under  Assumptions \ref{hyp_vector:indep_entries}, \ref{hyp:density} and \ref{hyp_vector:moments}  for $m=12$, the sequence $(\dist(B^n))_{n\ge 1}$ is $C$-tight, i.e.  is tight and has only  $C([0,1]^2)$-supported accumulation points.\end{lem}

 So let us now prove  Theorem  \ref{main_th0} under  Assumptions \ref{hyp_vector:indep_entries} and \ref{hyp_vector:moments2}. Note that it suffices to prove that the sequence $({distribution}(B^n))_{n\ge 1}$  has a unique possible accumulation  point supported by $C([0,1]^2)$ and that this accumulation  point  is the distribution of a centered Gaussian process which depends on the distributions of the $x_{i,j}$'s only through  $\lim_{n\to\infty}\E[|x_{1,2}|^4]$ (and actually does depend on $\lim_{n\to\infty} \E[|x_{1,2}|^4]$). Indeed, in this case, by Theorem 1.1 of \cite{cat-alain10}, where the case of GOE and GUE matrices is treated, this limit distribution is the bivariate Brownian bridge \ssi $\lim_{n\to\infty}\E[|x_{1,2}|^4]$ is the same as for $X_n$ a GOE or GUE matrix, i.e.  equal to $4-\beta$.

The following proposition is the key of the proof, since it allows to transfer the problem from the eigenvectors to some more accessible objects: the weighted spectral distributions of $X_n$.  

\begin{propo}\label{rep_4411_17h58} Let $e_i$ denote the $i$th vector of the canonical basis.  To prove Theorem  \ref{main_th0}, it suffices to prove that each   finite-dimentional marginal  distribution of the process $$\lf(\sum_{1\le i\le ns}\lf(e_i^*X_n^{k}e_i-\ff{n}\Tr X_n^{k}\ri)\ri)_{s\in [0,1], k\ge 1}$$  converges to  a centered Gaussian measure  and that the covariance of the limit process  depends on the distributions of the $x_{i,j}$'s only through $\lim_{n\to\infty} \E[|x_{1,2}|^4]$ and actually does depend on $\lim_{n\to\infty} \E[|x_{1,2}|^4]$.
\end{propo}

\begin{pr}Let $\mu_{X_n}:=\ff{n}\sum_{i=1}^n\delta_{\lam_i}$ be the empirical spectral law of $X_n$ and $F_{\mu_{X_n}}(u):=\ff{n}\sh\{i\ste \lam_i\le u\}$ be its cumulative distribution function. 
It is a well known result \cite[Th. 2.5]{bai-silver-book} that $F_{\mu_{X_n}}$ converges in probability, as $n$ tends to infinity, to the cumulative distribution function $F_{\scl}$ of the semicircle law in the space $D_c(\R, [0,1])$ (see Section \ref{sec:espaces+topo} for the definitions of the functional spaces and their topologies). 
Moreover,   the map \beq  D_0([0,1]^2) \times D_c(\R,[0,1]) & \longrightarrow & D_c([0,1]\times \R)\\
((G_{s,t})_{s,t\in [0,1]}, (g(u))_{u\in \R})& \longmapsto & (G_{s,g(u)})_{(s,u)\in [0,1]\times \R}
\eeq
is continuous at any pair of continuous functions. Hence for any continuous   process $(B_{s,t})_{s,t\in [0,1]}$ whose distribution  is an accumulation point of the sequence $(distribution(B^n))_{n\ge 1}$ for the Skorokhod topology in $D([0,1]^2)$,      the process $$(B^n_{s,F_{\mu_{X_n}}(u)})_{(s,u)\in [0,1]\times \R}$$ converges in distribution (up to the extraction of a subsequence) to the   process $$(B_{s,F_{\scl}(u)})_{(s,u)\in [0,1]\times \R}$$ (for the reader who is not used to weak convergence of \pro measures, this assertion relies on two results which can be found in \cite{billingsley}:  Theorem 4.4 and Corollary 1 of Theorem 5.1 in Chapter 1).
Now, note that $F_{\scl}:\R\to[0,1]$ admits a right inverse, so the distribution of the process $(B_{s,t})_{s,t\in [0,1]}$ is entirely determined by the one of the process $(B_{s,F_{\scl}(u)})_{(s,u)\in [0,1]\times \R}.$

As a consequence, to prove Theorem  \ref{main_th0}, it suffices to prove that the sequence     \be\la{30511.16h24}\lf(\dist((B^n_{s,F_{\mu_{X_n}}(u)})_{(s,u)\in [0,1]\times \R})\ri)_{n\ge 1}\ee has  a unique possible $C_c([0,1]\times \R)$-supported accumulation point as $n\to\infty$.

Now, note that any $f\in C_c([0,1]\times \R)$ is entirely determined by the collection of real numbers $(\int_{u\in \R} u^kf(s,u)\ud u)_{s\in [0,1], k\ge 0}$. Let us prove that this fact remains true in the ``distribution sense'' in the case where $f$ is random. More precisely, let us prove the following lemma.      

\begin{lem}\la{moments_pro_sto}
Let $f$ be a random variable taking values in $C_c([0,1]\times \R)$ \st with \pro one, $f(s,u)=0$ when $|u|>2$.  Then the distribution of $f$ is entirely determined by the  finite dimensional marginals of the process \be\label{3051116h00}
\lf(\int_{u\in \R} u^kf(s,u)\ud u\ri)_{s\in [0,1], k\ge 0}.\ee Moreover, in the case where the  finite dimensional marginals of the process of \eqref{3051116h00} are Gaussian and centered, then so are the ones of $f$.\end{lem}

\begin{pr}Let us fix $(s,u_0)\in [0,1]\times [-2,2]$ and let, for each $p\ge 1$,  $(P_{p,q})_{q\ge 1}$ be  a sequence of polynomials that is uniformly bounded on $[-3,3]$ and that converges pointwise  to $\one_{[u_0,u_0+1/p]}$ on $[-3,3]$. Then one has, with \pro one,  $$f(s,u_0)=\lim_{p\to\infty}p\int_{u_0}^{u_0+1/p}f(s,u)\ud u =\lim_{p\to\infty}\lim_{q\to\infty} p\int_{u\in \R}P_{p,q}(u)f(s,u)\ud u.$$This proves the lemma, because any almost sure limit   of a sequence of variables belonging to a space of  centered Gaussian variables is Gaussian and centered.
\end{pr}
 
 Since the fourth moment of the entries of $X_n$ is bounded, by e.g. \cite[Th. 5.1]{bai-silver-book}, we know that the extreme eigenvalues of $X_n$ converge to $-2$ and $2$. As a consequence, for any random variable $f$ taking values in $C_c([0,1]\times \R)$ \st the distribution of $f$ is a limit point of the sequence of  \eqref{30511.16h24}, we know that with \pro one, $f(s,u)=0$ when $|u|>2$.  
 
As a consequence, it follows from the previous lemma and from what precedes that to prove Theorem  \ref{main_th0}, it suffices to prove that each finite dimensional marginal  distribution of the process  $$\lf(\int_{u\in \R} u^kB^n_{s,F_{\mu_{X_n}}(u)}\ud u\ri)_{s\in [0,1],k\ge 0}$$ converges to   a centered Gaussian measure     and that the covariance of the limit process  depends on the distributions of the $x_{i,j}$'s only through $\lim_{n\to\infty} \E[|x_{1,2}|^4]$ and actually does depend on $\lim_{n\to\infty} \E[|x_{1,2}|^4]$.

But  for all $(s,u)\in [0,1]\times \R$, $$B^n_{s, F_{\mu_{X_n}}(u)}=\sqrt{\f{\beta}{2}}\sum_{1\le i\le ns}\sum_{\substack{1\le j\le n\\ \stm \lam_j\le u}}(|u_{i,j}|^2-1/n),$$hence 
 \be\label{441118h}B^n_{s, F_{\mu_{X_n}}(u)}=\sqrt{\f{\beta}{2}}\sum_{1\le i\le ns}F_{\mu_{X_n,e_i}-\mu_{X_n}}(u),\ee
where: \begin{itemize}\item $\mu_{X_n}$ is still the empirical spectral law of $X_n$,\item  $\mu_{X_n,e_i}$ is the weighted spectral law of $X_n$, defined by  $\mu_{X_n,e_i}:=\sum_{j=1}^n |u_{i,j}|^2\delta_{\lam_j}$,
\item $F_{\mu_{X_n,e_i}-\mu_{X_n}}$ is the cumulative distribution function of the null-mass signed measure ${\mu_{X_n,e_i}-\mu_{X_n}}$.\end{itemize}
 The following lemma  will allow to conclude the proof of Proposition \ref{rep_4411_17h58}. 
 \begin{lem}\label{moments-mesuredemassenulle}Let $\mu$ be a compactly supported null-mass signed measure and set $F_\mu(u):=\mu((-\infty, u])$. Then for all $k\ge 0$, $$\int_{u\in \R} u^kF_\mu(u)\ud u=-\int_{x\in \R} \f{x^{k+1}}{k+1}\ud\mu(x).$$\end{lem}
 \begin{pr}
 Let $a<b$ be such that the support of $\mu$ is contained in the open interval $(a,b)$. $F_\mu$ is null out of  $(a,b)$ and satisfies $F_\mu(u)=-\mu((u,b))$, so by Fubini's Theorem,    $$ \int_{u\in \R}\!\!\!\!\!\!\!\!  u^k F_\mu(u)\ud u=   \int_{x=a}^b\int_{u=a}^{x}\!\!\!\!\!\!\!\!   -u^k \ud u\ud \mu(x) 
 =\int_{x=a}^b \f{-x^{k+1}}{k+1}\ud\mu(x)+\f{a^{k+1}}{k+1}\mu((a,b))=\int_{x=a}^b \f{-x^{k+1}}{k+1}\ud\mu(x).$$
 \end{pr}
 
 It follows from this lemma and from  \eqref{441118h} that for all $s\in [0,1]$, $k\ge 0$, $$\int_{u\in \R} u^kB^n_{s,F_{\mu_{X_n}}(u)}\ud u=\f{-1}{k+1}\sqrt{\f{\beta}{2}}\sum_{1\le i\le ns}\lf(e_i^*X_n^{k+1}e_i-\ff{n}\Tr X_n^{k+1}\ri), 
 $$which proves Proposition \ref{rep_4411_17h58}.\end{pr}

Now, Theorem    \ref{main_th0} is a direct consequence of    the following proposition, whose proof is postponed below to Section \re{proofprop641114h}.
 
\begin{propo}\label{conv_moments}Under  Assumptions \ref{hyp_vector:indep_entries} and \ref{hyp_vector:moments2}, each   finite-dimentional marginal  distribution of the process $$\lf(\sum_{1\le i\le ns} (e^*_iX_n^ke_i-\ff{n}\Tr(X_n^k))\ri)_{s\in [0,1], k\ge 1}$$ converges to   a centered  Gaussian measure. The covariance of the limit distribution,  denoted by  $$\lf(\Co_{s_1,s_2}(k_1,k_2)\ri)_{s_1,s_2\in [0,1], k_1,k_2\ge 1},$$     depends on the distributions of the $x_{i,j}$'s only through  $\lim_{n\to\infty}\E[|x_{1,2}|^4]$.  Moreover, we have$$\Co_{s_1,s_2}(2,2)=(\lim_{n\to\infty}\E[|x_{1,2}|^4]-1)(\min\{s_1,s_2\}-s_1s_2).$$\end{propo}

\section{Proofs of the technical results}\la{sectechresults}
\subsection{Functional spaces and associated topologies}\label{sec:espaces+topo}
In this paper, we       use several functional spaces:\begin{itemize}\item
 $C([0,1]^2)$ (resp. $C_c([0,1]\times \R)$) is the set of continuous   functions on $[0,1]^2$ (resp. {\it compactly supported} continuous function on  $[0,1]\times \R$), endowed with the uniform convergence topology.\\
 \item $D_c(\R, [0,1])$ is the set of {\it compactly supported} c\`adl\`ag functions on $\R$ taking values in $[0,1]$, endowed with the topology defined by the fact that $f_n\lto f$ \ssi the bounds of the support of $f_n$ tend to the ones of the support of $f$ and for all $M>0$, after restriction to  $[-M,M]$, $f_n\lto f$ (the topology of $D([-M,M])$ being deduced from the one of $D([0,1])$ defined in   \cite{billingsley}).\\
 \item $D([0,1]^2)$ (resp. $D_c([0,1]\times \R)$)
 is the set   of functions $f:[0,1]^2\to \R$ (resp. of {\it compactly supported} functions  $f:[0,1]\times \R\to \R$) admitting limits in all ``orthants", more precisely \st for each $s_0,t_0$, for each pair of symbols $\diamond, \diamond'\in\{<,\ge\}$, 
 $$
 \lim_{\substack{s\diamond s_0\\ t\diamond' t_0}}f(s,t)$$ exists, and is equal to $f(s_0,t_0)$ if both $\diamond$ and $\diamond'$ are $\ge $. The space    $D([0,1]^2)$ is endowed with the Skorokhod topology
    defined in \cite{bickel} and the space $D_c([0,1]\times \R)$ is endowed with the topology defined by $f_n\lto f$ \ssi   for all $M>0$, after restriction to  $[0,1]\times [-M,M]$, $f_n\lto f$ (the topology of $D([0,1]\times [-M,M])$ being deduced from the one of $D([0,1]^2)$).\\
    \item $D_0([0,1]^2)$ is the set of functions in $D([0,1]^2)$ vanishing at the border of $[0,1]^2$, endowed with the induced topology.
\end{itemize}

\subsection{Proof of Proposition \ref{conv_moments}}\la{proofprop641114h}
Note that  by invariance of the law of $X_n$ under conjugation by any permutation matrix, the expectation of the random  law $\mu_{X_n,e_i}$ does not depend on $i $. So   for all $s\in [0,1]$, $ k\ge 1$,   \be\label{decB25211} \sum_{1\le i\le ns}(e^*_iX_n^ke_i-\ff{n}\Tr(X_n^k))=\sum_{1\le i\le ns}(e^*_iX_n^ke_i-\E[e^*_iX_n^ke_i])-\f{\lfloor ns\rfloor}{n}\sum_{1\le i\le n}(e^*_iX_n^ke_i-\E[e^*_iX_n^ke_i]).\ee
 
 Hence we are led to study the limit, as $n\to\infty$, of the finite-dimentional marginal distributions of the process \be\label{23.02.11.4a}\lf(\sum_{1\le i\le ns}(e^*_iX_n^ke_i-\E[e^*_iX_n^ke_i])\ri)_{s\in [0,1], k\ge 1}\ee  
 
The random variables of \eqref{23.02.11.4a} are going to be studied via their joint moments. So let us fix $p\ge 1$, $s_1, \ldots, s_p\in [0,1]$ and $k_1, \ldots, k_p\ge 1$. We shall study the limit, as $n$ tends to infinity, of \be\label{24211.moment}\E[\prod_{\ell=1}^p\sum_{1\le i\le ns_\ell} (e^*_iX_n^{k_\ell}e_i-\E[e^*_iX_n^{k_\ell}e_i])].\ee 
We introduce the set \be\label{presE}\mc{E}:=\{\ovl{0}^1, \ldots\ldots,\ovl{k_1}^1\}\cup \cdots \cdots\cdots \cdots\cup \{\ovl{0}^p, \ldots\ldots,\ovl{k_p}^p\},\ee where the sets $\{\ovl{0}^1, \ovl{1}^1, \ldots\}$, $\{\ovl{0}^2, \ovl{1}^2, \ldots\}$, \ldots\ldots,  $\{\ovl{0}^p, \ovl{1}^p, \ldots\}$  are $p$ disjoint copies of the set of nonnegative integers. The set $\mc{E}$ is ordered as presented in \eqref{presE}. For each partition $\pi$ of $\mc{E}$, for each $x\in \mc{E}$, we denote by $\pi(x)$   the  index of the class
of $x$,  after having ordered the classes according to the order of their first element (for example, $\pi(\ovl{1}^1)=1$;  $\pi(\ovl{2}^1)=1$ if $\ovl{1}^1\stackrel{\pi}{\sim}\ovl{2}^1$ and  $\pi(\ovl{2}^1)=2$ if $\ovl{1}^1\stackrel{\pi}{\nsim}\ovl{2}^1$).

By Assumption \ref{hyp_vector:indep_entries}, the expectation of \eqref{24211.moment}  can be   expanded and expressed as a sum on the set $\Part(\mc{E})$ of partitions of  the set $\mc{E}$ introduced above. 
We get  \beqy\label{24211.meddle}&\E[\prod_{\ell=1}^p\sum_{1\le i\le ns_\ell} (e^*_iX_n^{k_\ell}e_i-\E[e^*_iX_n^{k_\ell}e_i])]=&\\ \nonumber&n^{-\f{k_1+\cdots+k_p}{2}}\underset{\pi\in \Part(\mc{E})}{\sum}A(n,\pi)\E\lf[\prod_{\ell=1}^p\lf(x_{\pi(\ovl{0}^\ell),\pi(\ovl{1}^\ell) }\cdots x_{\pi(\ovl{k_\ell-1}^\ell),\pi(\ovl{k_\ell}^\ell) } -\E[x_{\pi(\ovl{0}^\ell),\pi(\ovl{1}^\ell) }\cdots x_{\pi(\ovl{k_\ell-1}^\ell),\pi(\ovl{k_\ell}^\ell) }]\ri)\ri],&
\eeqy
where for each $\pi\in \Part(\mc{E})$, $A(n,\pi)$ is the number of families of indices of   $$(i_{\ovl{0}^1}, \ldots,i_{\ovl{k_1}^1},i_{\ovl{0}^2}, \ldots,i_{\ovl{k_2}^2},\ldots\ldots\ldots,i_{\ovl{0}^p}, \ldots,i_{\ovl{k_p}^p})\in \ensn^{\mc{E}}$$ whose level sets partition
is $\pi$ and who satisfy, for each $\ell=1, \ldots,p$,   \be\label{54119h32}1\le  i_{\ovl{0}^\ell}=i_{\ovl{k_\ell}^\ell}\le ns_\ell.\ee

For any $\pi\in \Part(\mc{E})$, let us define $G_\pi$ to be the graph with vertex set the set $\{\pi(x)\ste x\in \mc{E}\}$ and edge set $$E_\pi:=\{\;\{\pi(\ovl{m-1}^\ell), \pi(\ovl{m}^\ell)\}\;\ste\; 1\le \ell\le p, \;m\in\{1, \ldots,k_\ell\}\;\}.$$  

For the term associated to $\pi$ in  \eqref{24211.meddle} to be non zero, we need to have:
\begin{itemize}
\item[(i)] $\pi(\ovl{0}^\ell)=\pi(\ovl{k_\ell}^\ell)$ for each $\ell=1, \ldots,p$, 
\item[(ii)] each edge of $G_\pi$ is visited at least twice by the union of the $p$ paths $(\pi(\ovl{0}^\ell), \ldots\ldots, \pi(\ovl{k_\ell}^\ell))$ ($\ell=1, \ldots,p$),
\item[(iii)] for each $\ell=1, \ldots, p$, the exists $\ell'\ne \ell$ \st at least one edge of $G_\pi$ is visited by both paths $(\pi(\ovl{0}^\ell), \ldots\ldots, \pi(\ovl{k_\ell}^\ell))$ and $(\pi(\ovl{0}^{\ell'}), \ldots\ldots, \pi(\ovl{k_{\ell'}}^{\ell'}))$. 
\end{itemize} 
Indeed, (i) is due to \eqref{54119h32}, (ii) is due to the fact that the $x_{i,j}$'s are centered and (iii) is due to the fact the the $x_{i,j}$'s are independent. 

Let us define a function $s(\cdot)$ on the set $\mc{E}$ in the following way: for each $\ell=1, \ldots,p$ and each $m=0, \ldots, k_\ell$, set  $$s({\ovl{m}^\ell})=\begin{cases}s_\ell&\textrm{if $m=0$ or $m=k_\ell$,}\\
1&\textrm{otherwise,}
\end{cases}$$ and   \be\label{24211.meddle.da2.1}s_\pi:=\prod_{B\textrm{ bloc of }\pi}\min_{x\in B}s(x).\ee  Then one can easily see that, as $n\to\infty$,   \be\label{24211.meddle.da2}A(n, \pi)\sim s_\pi n^{|\pi|}.\ee Thus for $\pi$ to have a non zero asymptotic contribution to \eqref{24211.meddle}, we need the following  condition, in addition to (i), (ii) and (iii):  \begin{itemize}
\item[(iv)]  $ \f{k_1+\cdots+k_p}{2}\le |\pi|$.\end{itemize}

The following lemma is a generalization of \cite[Lem. 2.1.34]{alice-greg-ofer}. Its proof goes along the same lines as the proof of the former. 
\begin{lem}Let $\pi\in \Part(\mc{E})$ satisfy (i), (ii), (iii). The the number $c_\pi$ of connected components of $G_\pi$ is \st $c_\pi\le p/2$ and \bes\label{maj_pi}|\pi|\le c_\pi-\f{p}{2}+\f{k_1+\cdots+k_p}{2}.\ees\end{lem}

As a consequence, if $\pi$ also satisfies (iv), we have \begin{enumerate}
\item[(a)] $c_\pi=p/2$,
\item[(b)] $p$ is even,
\item[(c)] $ |\pi |=\f{k_1+\cdots+k_p}{2} $ (so that $k_1+\cdots+k_p$ is also even).
\end{enumerate}
Note also that by (ii), we also have \begin{enumerate}
\item[(d)] $|E_\pi|\le \f{k_1+\cdots+k_p}{2}$. \end{enumerate}

To sum up, by \eqref{24211.meddle} and \eqref{24211.meddle.da2}, we have 
\beqy\label{25211.15h43}&\lim_{n\to\infty}\E[\prod_{\ell=1}^p\sum_{1\le i\le ns_\ell} (e^*_iX_n^{k_\ell}e_i-\E[e^*_iX_n^{k_\ell}e_i])]=&\\ \nonumber& \sum_{\pi}s_\pi\lim_{n\to\infty}\E\lf[\prod_{\ell=1}^p\lf(x_{\pi(\ovl{0}^\ell),\pi(\ovl{1}^\ell) }\cdots x_{\pi(\ovl{k_\ell-1}^\ell),\pi(\ovl{k_\ell}^\ell) } -\E[x_{\pi(\ovl{0}^\ell),\pi(\ovl{1}^\ell) }\cdots x_{\pi(\ovl{k_\ell-1}^\ell),\pi(\ovl{k_\ell}^\ell) }]\ri)\ri],&
\eeqy where the sum is taken over the partitions $\pi$ of $\mc{E}$ which satisfy (i), (ii), (iii) and (iv) above, and such partitions also do  satisfy (a), (b), (c) and (d) above.  

\noindent{\bf Case where $p$ is odd:} By (b), we know that when $p$ is odd, there is no partition $\pi$ satisfying the above conditions. So by \eqref{25211.15h43}, $$\lim_{n\to\infty}\E[\prod_{\ell=1}^p\sum_{1\le i\le ns_\ell} (e^*_iX_n^{k_\ell}e_i-\E[e^*_iX_n^{k_\ell}e_i])]= 0.$$

\noindent{\bf Case where $p=2$:} In this case,  by (a), for each partition $\pi$ satisfying (i), (ii), (iii) and (iv), $G_\pi$ is connected, so that $|\pi|-1\le |E_\pi|$. Hence by (c) and (d),  $|E_\pi|$ is either equal to $\f{k_1+ k_2}{2}$ or to $\f{k_1+k_2}{2}-1$.

In the case where $| E_\pi|=\f{k_1+k_2}{2}-1$, the graph $G_\pi$ has exactly one more vertex than edges, hence is a tree. As a consequence, the paths $(\pi(\ovl{0}^1), \ldots\ldots, \pi(\ovl{k_1}^1))$ and
$ (\pi(\ovl{0}^2), \ldots\ldots, \pi(\ovl{k_2}^2))$,
which have the same beginning and ending vertices, have the property to visit an even number of times each edge they visit. By an obvious cardinality argument, only one edge is visited more than twice, and it is visited four times (twice in each sense). The other edges are visited once in each sense. It follows that the expectation associated to $\pi$ in  \eqref{25211.15h43} is equal to $\E[|x_{1,2}|^4]-1$.

In the case where $| E_\pi|=\f{k_1+k_2}{2}$,  by a  cardinality argument again, we see that each edge of $G_\pi$ is visited exactly twice (in fact, the configuration is the one described in  \cite[Sect. 2.1.7]{alice-greg-ofer}, where $G_\pi$ is a {\it  bracelet}). It follows that the expectation associated to $\pi$ in  \eqref{25211.15h43} is equal to $1$.

As a consequence, as $n$ tends to infinity, $$\E \lf[\sum_{1\le i\le ns_1} (e^*_iX_n^{k_1}e_i-\E[e^*_iX_n^{k_1}e_i])\times  \sum_{1\le i\le ns_2} (e^*_iX_n^{k_2}e_i-\E[e^*_iX_n^{k_2}e_i])\ri]$$ converges to a number that we shall denote by  \be\label{defcov}\Co^{\textrm{centered}}_{s_1,s_2}(k_1,k_2)\ee and  which  depends on the distributions of the $x_{i,j}$'s only through $\lim_{n\to\infty}\E[|x_{1,2}|^4]$. 
\\

\noindent{\bf Case where $p$ is $>2$ and even:} By (a) above, for each partition $\pi$ satisfying (i), (ii), (iii) and (iv),  $G_\pi$ has exactly $p/2$ connected components. By (iii), each of them contains the support of exactly two of the $p$ paths \be\label{25211.1}   \qquad \qquad\qquad\qquad\qquad(\pi(\ovl{0}^\ell), \ldots\ldots, \pi(\ovl{k_\ell}^\ell))  \qquad\qquad \qquad\qquad\qquad \textrm{ ($\ell=1, \ldots,p$).}\ee Let us define $\si_\pi$ to be the {\it matching} (i.e. a permutation all of whose cycles have length two) of $\{1, \ldots, p\}$ \st for all $\ell=1, \ldots, p$, the paths with indices $\ell$ and $\si_\pi(\ell)$ of the collection  \eqref{25211.1} above are supported by the same connected component of  $G_\pi$. 

We shall now partition the sum of  \eqref{25211.15h43}   according to the value of the matching $\si_\pi$ defined by $\pi$. We get 
\beq&\lim_{n\to\infty}\E[\prod_{\ell=1}^p\sum_{1\le i\le ns_\ell} (e^*_iX_n^{k_\ell}e_i-\E[e^*_iX_n^{k_\ell}e_i])]=&\\ \nonumber&\sum_\si\sum_\pi s_\pi\lim_{n\to\infty}\E\lf[\prod_{\ell=1}^p\lf(x_{\pi(\ovl{0}^\ell),\pi(\ovl{1}^\ell) }\cdots x_{\pi(\ovl{k_\ell-1}^\ell),\pi(\ovl{k_\ell}^\ell) } -\E[x_{\pi(\ovl{0}^\ell),\pi(\ovl{1}^\ell) }\cdots x_{\pi(\ovl{k_\ell-1}^\ell),\pi(\ovl{k_\ell}^\ell) }]\ri)\ri],&
\eeq
where the first sum is over the matchings $\si$ of $\{1, \ldots, p\}$  and the second sum is over the partitions $\pi$ satisfying (i), (ii), (iii) and (iv) \st $\si_\pi=\si$. 

Note that for each matching $\si$ of $\{1, \ldots, p\}$, the set of   partitions $\pi$ of $\mc{E}$ \st $\si_\pi=\si$ can be identified with the cartesian product, indexed by the set  of cycles $\{\ell,\ell'\}$ of $\si$, of the set of partitions $\tau$ of the set $$\mc{E}_{\ell,\ell'}:=\{\ovl{0}^\ell, \ldots\ldots,\ovl{k_\ell}^\ell\} \cup \{\ovl{0}^{\ell'}, \ldots\ldots,\ovl{k_{\ell'}}^{\ell'}\}\qquad\textrm{ (subset of $\mc{E}$)}$$
satisfying the following conditions  \begin{itemize}
\item[(i.2)] $\tau(\ovl{0}^\ell)=\tau(\ovl{k_\ell}^\ell)$  and $\tau(\ovl{0}^{\ell'})=\tau(\ovl{k_{\ell'}}^{\ell'})$, 
\item[(ii.2)] each edge of the graph $G_\tau$  is visited at least twice by the union of the 2 paths $(\tau(\ovl{0}^\ell), \ldots\ldots, \tau(\ovl{k_\ell}^\ell))$ and $(\tau(\ovl{0}^{\ell'}), \ldots\ldots,\tau(\ovl{k_{\ell'}}^{\ell'}))$,
\item[(iii.2)]  at least one edge of $G_\tau$ is visited by both previous paths, 
\item[(iv.2)]  $ \f{k_\ell+ k_{\ell'}}{2}\le |\tau|$.\end{itemize}
Moreover, by independence of the random variables $x_{i,j}$'s, the expectation $$\E\lf[\prod_{\ell=1}^p\lf(x_{\pi(\ovl{0}^\ell),\pi(\ovl{1}^\ell) }\cdots x_{\pi(\ovl{k_\ell-1}^\ell),\pi(\ovl{k_\ell}^\ell) } -\E[x_{\pi(\ovl{0}^\ell),\pi(\ovl{1}^\ell) }\cdots x_{\pi(\ovl{k_\ell-1}^\ell),\pi(\ovl{k_\ell}^\ell) }]\ri)\ri]$$ factorizes along the connected components of $G_\pi$. The factor $s_\pi$, defined in  \eqref{24211.meddle.da2.1},   also factorizes along the connected components of $G_\pi$. It follows that we have \beq&\lim_{n\to\infty}\E[\prod_{\ell=1}^p\sum_{1\le i\le ns_\ell} (e^*_iX_n^{k_\ell}e_i-\E[e^*_iX_n^{k_\ell}e_i])]=&\\ \nonumber& {\sum_\si} \prod_{\{\ell,\ell'\}}\lim_{n\to\infty}\E\lf[\sum_{1\le i\le ns_\ell} (e^*_iX_n^{k_\ell}e_i-\E[e^*_iX_n^{k_\ell}e_i])\times \sum_{1\le i\le ns_{\ell'}} (e^*_iX_n^{k_{\ell'}}e_i-\E[e^*_iX_n^{k_{\ell'}}e_i])\ri],&
\eeq
where the sum is over the matchings $\si$ of $\{1, \ldots, p\}$ and for each such $\si$, the product is over the cycles $\{\ell,\ell'\}$ of $\si$.

By the definition of $\Co^{\textrm{centered}}_{s_1,s_2}(k_1,k_2)$ in \eqref{defcov}, we get  $$\lim_{n\to\infty}\E[\prod_{\ell=1}^p\sum_{1\le i\le ns_\ell} (e^*_iX_n^{k_\ell}e_i-\E[e^*_iX_n^{k_\ell}e_i])]= \sum_{\si\textrm{ matching  }}\prod_{\textrm{ $\{\ell,\ell'\}$ cycle of $\si$}}\Co^{\textrm{centered}}_{s_\ell,s_{\ell'}}(k_\ell,k_{\ell'}).$$ By Wick's formula and Equation \eqref{decB25211}, we have proved the first part of Proposition \ref{conv_moments}. 
\\

\noindent{\bf Computation of $\Co_{s_1,s_2}(2,2)$:} We have, by the paragraph devoted to the case $p=2$ above, $$\Co^{\textrm{centered}}_{s_1,s_2}(2,2)=\sum_{\pi, \, G_\pi\textrm{ is a tree}}s_\pi\times(\lim_{n\to\infty}\E[|x_{1,2}|^4]-1)\;\;+\sum_{\pi, \, G_\pi\textrm{ is a bracelet}}s_\pi.$$
There is exactly one tree and zero bracelet with 2 vertices. We represent this tree in the following way: $\bullet-\circ$. There are two associated partitions $\pi$:\begin{itemize}
\item[-] the first one is defined by $$(\pi(\ovl{0}^1),\pi(\ovl{1}^1),\pi(\ovl{2}^1))=(\pi(\ovl{0}^2),\pi(\ovl{1}^2),\pi(\ovl{2}^2))=(\bullet, \circ, \bullet),$$and satisfies $s_\pi=\min\{s_1,s_2\}$,
\item[-] the second one is defined by $$(\pi(\ovl{0}^1),\pi(\ovl{1}^1),\pi(\ovl{2}^1))=(\bullet, \circ, \bullet) \;,\qquad   (\pi(\ovl{0}^2),\pi(\ovl{1}^2),\pi(\ovl{2}^2))=(\circ, \bullet, \circ),$$and satisfies $s_\pi=s_1s_2$.
\end{itemize}
As a consequence, $\Co^{\textrm{centered}}_{s_1,s_2}(2,2)=(\lim_{n\to\infty}\E[|x_{1,2}|^4]-1)(\min\{s_1,s_2\}+s_1s_2).$ 

From this formula and Equation \eqref{decB25211}, we easily deduce the formula of $\Co_{s_1,s_2}(2,2)$. It concludes the proof of Proposition \ref{conv_moments}.\hfill$\square$

\subsection{A preliminary result for the proofs of Proposition \ref{bound8411} and Lemma \ref{tightness_lemma}}
 Proposition \ref{conv_esp_1411} below will allow to prove Lemma \ref{tightness_lemma}.   Let us   now consider two independent Wigner matrices $X_n,X_n'$.    Let us introduce a modified version of Assumption \ref{hyp_vector:moments}, 	also depending on an integer $m\ge 2$:
  \begin{assum}\label{hyp_vector:momentsprime} For each $k\ge 1$, $\sup_n \E[|x_{1,1}|^{k}+|x_{1,2}|^{k}+|x'_{1,1}|^{k}+|x'_{1,2}|^{k}]<\infty$. 
  Moreover, there exists $\eps_0>0$ \st   for each $r,s\ge 0$,  \be\label{bound_moments31311pp4}r+s\le m-2\;\Longrightarrow\; \E[\Re(x_{1,1})^r\Im(x_{1,1})^s]-\E[\Re(x'_{1,1})^r\Im(x'_{1,1})^s]=O(n^{-\eps_0+1+\f{r+s-m}{2}})\ee and \be\label{bound_moments31311pq5}r+s\le m\;\Longrightarrow \;\E[\Re(x_{1,2})^r\Im(x_{1,2})^s]-\E[\Re(x'_{1,2})^r\Im(x'_{1,2})^s]=O(n^{-\eps_0+\f{r+s-m}{2}}).\ee
\end{assum}

 Let $U_n'=[u_{i,j}']_{i,j=1}^n$ be the associated eigenvectors matrix of a  matrix $X'_n$ (like $U_n=[u_{i,j}]_{i,j=1}^n$ for $X_n$).  
\begin{propo}\label{conv_esp_1411} We suppose that $X_n$, $X_n'$ both satisfy Assumptions \ref{hyp_vector:indep_entries} and   \ref{hyp:density}, and satisfy Assumption \ref{hyp_vector:momentsprime}.
Let us fix a positive integer $k$ and a polynomial function $G$ on $\C^k$. For each $n$, let us consider a collection $(i_1,p_1,q_1), \ldots, (i_k,p_k,q_k)\in \{1, \ldots, n\}^3$ (this collection might depend on $n$). Then for a certain constant $C$ independent of $n$ and of the choices of $(i_1,p_1,q_1), \ldots, (i_k,p_k,q_k)$ (but depending on $k$ and on $G$)
\be\label{ineqTaoVurate}|\;\E[G(n\,u_{p_1,i_1}\ovl{u}_{q_1,i_1}, \ldots,n\,u_{p_k,i_k}\ovl{u}_{q_k,i_k})]-\E[G(n\,u'_{p_1,i_1}\ovl{u}'_{q_1,i_1}, \ldots,n\,u'_{p_k,i_k}\ovl{u}'_{q_k,i_k})]\;|\le Cn^{2-\f{m}{2}}.\ee
\end{propo}

\begin{pr}
The proof follows the ideas developed by Terence Tao and Van Vu in recent papers, such as \cite{Tao-Vu_0906.0510,Tao-Vu_1103.2801v1}. More precisely, it makes an intensive use of the strategy and of some estimations given   in \cite{Tao-Vu_1103.2801v1}. 

An event $E=E^{(n)}$ depending on the parameter $n$ will be said to be true with {\it overwhelming probability} if for all $c>0$, $1-\Pro(E)=O(n^{-c})$. One can neglect any  such event in the proof of \eqref{ineqTaoVurate} because $G$ has polynomial growth and the entries of unitary matrices are bounded by $1$.

  
 By Assumption \ref{hyp_vector:moments} and Chebyshev's inequality, that for each $k$, as $t\to\infty$, $$\Pro(  |x_{i,j}|\ge t)=O(t^{-k})\qquad\textrm{  (uniformly in $i,j$ and   $n$),}$$so that by the union bound,  $\Pro(\max_{1\le i,j\le n}|x_{i,j}|\ge n^\eps)=O( n^{2-k\eps})$ for any $\eps,k$. Hence by what precedes, for a fixed $\eps_1$ (a parameter fixed later), one can suppose that for each $n$,   the $x_{i,j}$'s satisfy $|x_{i,j}|\le n^{\eps_1}$ for all $i,j$. Of course, one can also suppose that $|x'_{i,j}|\le n^{\eps_1}$ for all $i,j$.
 

 
 To prove \eqref{ineqTaoVurate}, it suffices to prove that we can replace the $x_{i,j}$'s by the $x'_{i,j}$'s one by one up to an error in the considered expectation which is $O(n^{1-\f{m}{2}})$ for each diagonal replacement and $O(n^{-\f{m}{2}})$ for each off-diagonal replacement (and, of course, that these  bounds on the errors are uniform in the $\f{n(n+1)}{2}$ replacements).  
  
 So let us fix $1\le p\le q\le n$ and define, for each $z\in \C $ (or $z\in \R$ if $p=q$),   $A_{p,q}(z)$ be the matrix obtained by ``mixing" the (rescaled) matrices $X_n$ and $X_n'$ in the following way:
\begin{itemize}
\item for each $1\le i\le j\le n$ \st $(i,j)$ is below $(p,q)$ for the lexicographic order, the $i,j$ entry of $A_{p,q}(z)$ and its symmetric one are the corresponding ones in the entries of $nX_n'$,
\item the $p,q$ entry of $A_{p,q}(z)$ and its symmetric one are $z$ and $\ovl{z}$,
\item for each $1\le i\le j\le n$ \st $(i,j)$ is above $(p,q)$ for the lexicographic order, the $i,j$ entry of $A_{p,q}(z)$ and its symmetric one are the corresponding ones in the entries of $nX_n$.\end{itemize}
 We set $U(z)=[u_{i,j}(z)]_{i,j=1}^n$ to be the eigenvectors matrix of $A_{p,q}(z)$ and $$F(z):=G(n\,u_{p_1,i_1}(z)\ovl{u}_{q_1,i_1}(z), \ldots,n\,u_{p_k,i_k}(z)\ovl{u}_{q_k,i_k}(z)).$$  
 Let $\E[\;\cdot\;|\;A_{p,q}(0)]$ denote the conditional expectation with respect to the $\si$-algebra generated by the entries of $A_{p,q}(0)$. We have to prove  that with overwhelming probability, we have, uniformly on the \pro space and on $p,q$, \be\label{ineqTaoVurateprime}\E[F(\sqrt{n}x_{p,q})\;|\;A_{p,q}(0)]=\E[F(\sqrt{n}x'_{p,q})\;|\;A_{p,q}(0)]+O(n^{\one_{p=q}-\f{m}{2}}).\ee
  To prove it, we shall use the Taylor expansion of $F$ around zero. Let us for example treat the case $p\ne q$. 
  For $|x|\le n^{\eps_1}$, we have 
  $$
  F(\sqrt{n}x)=\sum_{\substack{  r+s\le m}}\f{n^{\f{r+s}{2}}}{r!s!} \partial_{\Re(z)}^{r}\partial_{\Im(z)}^s F(0)\Re(x)^s \Im(x)^s+O\lf(n^{(m+1)({\eps_1+1/2})}\sup_{\substack{r+s=m+1\\ |z|\le n^{\eps_1+1/2}}}|\partial_{\Re(z)}^{r}\partial_{\Im(z)}^sF(z)|\ri).
  $$
 Moreover,  by \cite[Propo. 20]{Tao-Vu_1103.2801v1}, we know that with overwhelming probability, $A_{p,q}(0)$ is a {\it good configuration} (in the sense of \cite[Def. 18]{Tao-Vu_1103.2801v1}). By Equation (27) of \cite[Lem. 24]{Tao-Vu_1103.2801v1}, it implies that with overwhelming probability, for all $l=1, \ldots, k$, for all $r,s$ \st $r+s\le 10$, \bes\label{borne_der_31311}\sup_{|z|\le n^{1/2+\eps_1}} |\partial_{\Re(z)}^{r}\partial_{\Im(z)}^s(nu_{p_l,i_l}(z)\ovl{u}_{q_l,i_l}(z)|)=
 O(n^{\eps_1-(r+s)})
 \ees
  (in the statement of \cite[Lem. 24]{Tao-Vu_1103.2801v1}, $m$ must be bounded above by $10$, but the bound 10 can be replaced by any finite one). 
 The function $G$ is   polynomial, let $d$ be its degree. It
 follows 
 that with overwhelming probability, for all $r,s$, we have 
 $$  \sup_{|z|\le n^{1/2+\eps_1}} |\partial_{\Re(z)}^{r}\partial_{\Im(z)}^s F(z)|=O(n^{d\eps_1-(r+s)}).
 $$  As
 a consequence, by the bounds of \eqref{bound_moments31311pq5}, 
$$\E[F(\sqrt{n}x_{p,q})\;|\;A_{p,q}(0)]-\E[F(\sqrt{n}x'_{p,q})\;|\;A_{p,q}(0)]= \sum_{r+s\le m} O(n^{d\eps_1-\f{m}{2}-\eps_0})+O(n^{(m+1+d){\eps_1}-\f{m+1}{2}})$$
Now, choosing $\eps_1$ small enough so that $d\eps_1<\eps_0$ and $(m+1+d){\eps_1}<1/2$, we get \eqref{ineqTaoVurateprime}.
\end{pr}

\subsection{Proof of Proposition \ref{bound8411}}Let us first  prove that $\E[(n|u_{i,j}|^2-1)^4]$ is bounded independently of $i,j$ and $n$. To do so, let us 
    apply  Proposition \ref{conv_esp_1411} for $X_n'$ a GOE or GUE matrix and  $m=4$: it allows to reduce the problem to the case where the $u_{i,j}$'s are the entries of a Haar-distributed matrix. In this case, the distribution of $u_{i,j}$ does not depend on $i,j$ and it is well known that the moments of $\sqrt{n}u_{1,1}$ converge  to the ones of a standard real or complex Gaussian variable. We did not find any concise enough reference for this moments convergence,  but in the unitary case, Proposition 3.4 of \cite{FBDAOPinfdiv} allows to compute $\E[(n|u_{1,1}|^2-1)^4]$ and to verify its boundedness. To treat the orthogonal  case, the most direct way to compute $\E[(n|u_{1,1}|^2-1)^4]$ is to use the fact that $u_{1,1}^2$ has the same distribution as $Z_1^2/(Z_1^2+\cdots+Z_n^2)$ for $Z_i$ independent standard Gaussian variables  (see e.g. \cite[Lem. 2.1]{jang-variance-formula}) and then to use  \cite[Lem. 2.4]{jang-variance-formula}.
  
  Let us now prove that there is a constant $C$ independent of $n$ and of  $0\le s<s'\le 1$ and $0\le t<t'\le 1$ \st  \be\label{8411_21h43}\E[\{\sum_{\substack{ns<i\le ns'}}\sum_{{nt<j\le nt'}}(|u_{i,j}|^2-1/n)\}^2]\le C (s'-s)(t'-t).\ee
  Set $I=\{i=1, \ldots, n\ste ns<i\le ns'\}$ and $J=\{j=1, \ldots, n\ste nt<j\le nt'\}$. We have 
$$\E[
\{\sum_{i\in I,j\in J}(|u_{i,j}|^2-1/{n})\}^{2}]=n^{-2}\sum_{(i_\ell,j_\ell)_{\ell=1,2} \in (I\times J)^2}\E[\prod_{\ell=1,2}(n|u_{i_\ell,j_\ell}|^2-1)].$$
Let us now apply  Proposition \ref{conv_esp_1411} for $X_n'$ a GOE or GUE matrix and  $m=4$. Let $U'_n=[u'_{i,j}]_{i,j=1}^n$ be a Haar-distributed orthogonal or unitary  matrix. By  Proposition \ref{conv_esp_1411}, we have 
$$\E[
\{\sum_{i\in I,j\in J}(|u_{i,j}|^2-1/{n})\}^{2}]=n^{-2}O(n^{-2}(\sh I\times J)^2)+n^{-2}\sum_{(i_\ell,j_\ell)_{\ell=1,2} \in (I\times J)^2}\E[\prod_{\ell=1,2}(n|u'_{i_\ell,j_\ell}|^2-1)],$$ where the term $n^{-2}O(n^{-2}(\sh I\times J)^2)$ is uniform in $s,s',t,t'$, hence is bounded,
for a certain constant $C_1$, by $C_1(s'-s)^2(t'-t)^2\le C_1(s'-s)(t'-t)$. 
So  it suffices to prove the result for $U_n'$ instead of $U_n$.  
Since each $|u'_{i,j}|^2$ has expectation $1/n$,  $$
\E[\prod_{\ell=1}^2(n|u'_{i_\ell,j_\ell}|^2-1)]=
n^2\E[|u'_{i_1,j_1}|^2|u'_{i_2,j_2}|^2]-1.
$$ Then, to compute $\E[|u'_{i_1,j_1}|^2|u'_{i_2,j_2}|^2]$, one uses the {\it Weingarten calculus}, developed in \cite{collinsIMRN,collins-sniady06}. By the formulas of \cite[Sec. 6]{collins-sniady06} for the orthogonal group and  \cite[Sec. 5.2]{collinsIMRN} for the unitary group,   \begin{align*} &\op{Wg}^{\textrm{ortho}}(Id_{\{1,2\}})=\f{n+1}{n(n-1)(n+2)},&\quad& \op{Wg}^{\textrm{ortho}}((12))=\f{1}{n(n-1)(n+2)},\\
 &\op{Wg}^{\textrm{unit}}(Id_{\{1,2\}})=\f{1}{n^2-1},&\quad& \op{Wg}^{\textrm{unit}}((12))=\f{-1}{n(n^2-1)},\end{align*} (where $(12)$ denotes the transposition with support   $\{1,2\}$), one can easily verify that  $$n^2\E[|u'_{i_1,j_1}|^2|u'_{i_2,j_2}|^2]-1=\begin{cases}O(n^{-2})&\textrm{ if $i_1\ne i_2$, $j_1\ne j_2$,}\\
 O(n^{-1} ) &\textrm{ if ($i_1=i_2$, $j_1\ne j_2$) or ($i_1\ne i_2$, $j_1= j_2$),}\\O(1) &\textrm{ if $i_1=i_2$, $j_1= j_2$,}
 \end{cases}$$ which is enough to deduce \eqref{8411_21h43}.


\subsection{Proof of Lemma \ref{tightness_lemma}}\la{prooflem641114h} To prove Lemma  \ref{tightness_lemma}, we shall use the following proposition, which is the obvious multidimentional generalization of Proposition 3.26 of Chapter VI of \cite{jacod-shi} (note also that the result of \cite{jacod-shi} is devoted to pocesses on $[0,\infty)$ and not on $[0,1]$, so that the $N$ of \cite{jacod-shi} means nothing to us). 

For $f\in D([0,1]^2)$ and  $(s_0,t_0)\in [0,1]^2$, we define $\Delta_{s_0,t_0} f$ to be the {\it maximal jump} of $f$ at $(s_0,t_0)$, i.e. $$\Delta_{s_0,t_0} f:=\max_{\diamond, \diamond'\in\{<,\ge\}}\lf|f(s_0,t_0)-
 \lim_{\substack{s\diamond s_0\\ t\diamond' t_0}}f(s,t)\ri|.$$

\begin{propo}
If the sequence $(\dist(B^n))_{n\ge 1}$ is tight and satisfies \be\label{disccv0}\forall \eps >0,\qquad\qquad \Pro(\sup _{(s,t)\in [0,1]^2}\Delta_{s,t} B^n >\eps)\ninf 0,\ee then the sequence $(\dist(B^n))_{n\ge 1}$ is $C$-tight, i.e. is tight and can only have $C([0,1]^2)$-supported accumulation points.
\end{propo}

So to prove Lemma \ref{tightness_lemma}, let us first prove that the sequence $(distribution(B^n))_{n\ge 1}$ is tight.  Note that   the process $B^n$ vanishes at the border of $[0,1]^2$.  So   according to \cite[Th. 3]{bickel} and to Cauchy-Schwartz inequality,  
it suffices to prove that there is $C<\infty$    \st for $n$ large enough,  for all   $s<s', \, t<t'\in [0,1]$, \bes\label{majacc}\E[
\{\sum_{ns<i\le ns'}\sum_{nt< j\le nt'}(|u_{i,j}|^2-1/{n})\}^{4}]\le C(s'-s)^{2}(t'-t)^{2}.\ees
As in the proof of \eqref{8411_21h43} above, one can suppose that the $u_{i,j}$'s are the entries of a Haar-distributed matrix.
 But in this case, the job has already   been done in \cite{cat-alain10}: the unitary case is treated in  Section 3.4.1 (see specifically  Equation (3.25)) and the orthogonal case is treated, more elliptically,  in Section 4.5  (to   recover the details of the proof, join Equations (3.26), (4.5) and (4.17)).

Let us now prove \eqref{disccv0}. Note that $\sup _{(s,t)\in [0,1]^2}\Delta_{s,t} B^n =\max_{1\le i,j\le n} ||u_{i,j}|^2-1/n|$. As a consequence, by the union bound, it suffices to prove that for each $\eps>0$, there is $C<\infty$ independent of $i,j$ and $n$ \st for all $i,j$, \bes\label{majacc17h23}\Pro(||u_{i,j}|^2-1/n |>\eps)\le Cn^{-4},\ees which follows from 
 Chebyshev's inequality and Proposition \ref{bound8411}.\hfill$\square$


\begin{thebibliography}{10}
\bibitem{alice-greg-ofer}
G.~Anderson, A.~Guionnet, O.~Zeitouni
\newblock \emph{An Introduction to Random Matrices}.
\newblock  Cambridge studies in advanced mathematics, {118} (2009).

\bibitem{baipanJStatPhys2012} Z.D. Bai, G.M. Pan \emph{Limiting Behavior of Eigenvectors of Large Wigner Matrices},   J. Stat. Phys. 146(3) 519--549, 2012.

\bibitem{bai-silver-book} Z.~D.~Bai, J.~W.~Silverstein \emph{Spectral analysis of large dimensional random matrices}, Second Edition, Springer, New York, 2009.

\bibitem{FBDAOPinfdiv} {F. Benaych-Georges} \emph{Classical and free infinitely divisible distributions and random matrices}. {Annals of Probability}.  Vol. 33, no. 3 (2005) 1134--1170.

\bibitem{bickel} P.J. Bickel, M.J. Wichura \emph{Convergence criteria for multiparameter stochastic processes and some applications}, Ann. Math. Statist., 42(5):1656--1670, 1971.

\bibitem{billingsley} P. Billingsley  \emph{Convergence of Probability Measures}. New York, John 
Wiley and Sons (1968).  

\bibitem{chafai-anderson} D. Chafa\"{\i}   \emph{Anderson localization}. Blog note, available at \url{http://blog.djalil.chafai.net/2010/09/07/anderson-localization/}

\bibitem{meckessourav08}  S. Chatterjee, E.  Meckes  \emph{Multivariate normal approximation using exchangeable pairs},  ALEA 4, 257--283 (2008). 

\bibitem{collinsIMRN} B. Collins \emph{Moments and Cumulants of Polynomial Random Variables on Unitary Groups, the Itzykson-Zuber Integral, and Free Probability}. Int. Math. Res. Not., No. 17, 2003. 

\bibitem{collins-sniady06}  B. Collins, P. Sniady \emph{Integration with respect to the Haar measure on unitary, orthogonal and symplectic group}. Comm. Math. Phys., 264, 773--795, 2006.

 \bibitem{cat-alain10} C. Donati-Martin, A. Rouault \emph{Truncations of Haar unitary matrices, traces and bivariate Brownian bridge}, Random Matrices: Theory and Application (RMTA) Vol. 01, No. 01.

\bibitem{erdosknowlesband1} L. Erd\"os, A. Knowles \emph{Quantum Diffusion and Eigenfunction Delocalization in a Random Band Matrix Model}. Comm. Math. Phys. 303 (2011), no. 2, 509--554.

\bibitem{erdosknowlesband2} L. Erd\"os, A. Knowles \emph{Quantum Diffusion and Delocalization for Band Matrices with General Distribution}. Ann. Henri Poincaré 12 (2011), no. 7, 1227--1319.


\bibitem{ErdosSchleinYau1} L. Erd\"os, B. Schlein, H.-T. Yau
\emph{Semicircle law on short scales and delocalization of eigenvectors for Wigner random matrices}.
Ann. Probab. 37 (2009), no. 3, 815--852. 

\bibitem{ErdosSchleinYau2} L. Erd\"os, B. Schlein, H.-T. Yau
\emph{Local semicircle law and complete delocalization for Wigner random matrices}.
Comm. Math. Phys. 287 (2009), no. 2, 641-655. 


\bibitem{ErdosSchleinYau3} L. Erd\"os, B. Schlein, H.-T. Yau
\emph{Wegner estimate and level repulsion for Wigner random matrices}, 
Int. Math. Res. Not. IMRN 2010, no. 3, 436--479.

 \bibitem{jacod-shi} J. Jacod, A.N. Shiryaev \emph{Limit theorems for stochastic processes}, Berlin, Springer-Verlag (1987). 

\bibitem{jiang06}     T. Jiang  \emph{How many entries of a typical orthogonal matrix can be approximated by independent normals?} Ann. Probab. 34(4), 1497--1529. 2006.

\bibitem{jang-variance-formula}T. Jiang \emph{A variance formula related to a quantum conductance problem}. Physics
Letters A, vol. 373 (25), 2117--2121, 2009.

\bibitem{KnowlesYinEig} A. Knowles, J. Yin \emph{Eigenvector Distribution of Wigner Matrices}, arXiv:1102.0057, to appear in Prob. Theor. Rel. Fields.

\bibitem{schenker_band} J. Schenker	\emph{Eigenvector localization for random band matrices with power law band width}. 	Comm. Math. Phys. 290, (2009) 1065--1097.

\bibitem{silverstein-AOP-eigenvectors} J. W. Silverstein \emph{Weak convergence of random functions defined by the eigenvectors of sample covariance matrices}.  Ann. Probab. 18 (1990), no. 3, 1174--1194. 

\bibitem{s_sodin_band_matrices} S. Sodin \emph{The spectral edge of some random band matrices}. 
Ann. of Math. (2) 172 (2010), no. 3, 2223--2251. 

\bibitem{taovucircularlaw} T. Tao, V. Vu (with an appendix by  M. Krishnapur) \emph{Random matrices: Universality of ESDs and the circular law}, Ann. Probab. Volume 38, Number 5 (2010), 2023--2065.

\bibitem{Tao-Vu_0906.0510} T. Tao, V. Vu \emph{Random matrices: universality of local eigenvalue statistics},  Acta Mathematica, 	
 206 (2011), 127--204.

\bibitem{Tao-Vu_1103.2801v1} T. Tao, V. Vu \emph{Random matrices: Universal properties of eigenvectors}, Random matrices: Theory and Applications (RMTA), Vol. 01, No. 01. DOI: 10.1142/S2010326311500018
\end{thebibliography}
\end{document}